\documentclass[reqno]{amsart}
\usepackage{color}

\newtheorem{theorem}{Theorem}
\newtheorem{lemma}[theorem]{Lemma}

\theoremstyle{definition}

\theoremstyle{remark}

\numberwithin{equation}{section}

\newcommand{\intav}[1]{\mathchoice {\mathop{\vrule width 6pt height 3 pt depth  -2.5pt
\kern -8pt \intop}\nolimits_{\kern -6pt#1}} {\mathop{\vrule width
5pt height 3  pt depth -2.6pt \kern -6pt \intop}\nolimits_{#1}}
{\mathop{\vrule width 5pt height 3 pt depth -2.6pt \kern -6pt
\intop}\nolimits_{#1}} {\mathop{\vrule width 5pt height 3 pt depth
-2.6pt \kern -6pt \intop}\nolimits_{#1}}}

\newcommand{\intavl}[1]{\mathchoice {\mathop{\vrule width 6pt height 3 pt depth  -2.5pt
\kern -8pt \intop}\limits_{\kern -6pt#1}} {\mathop{\vrule width 5pt
height 3  pt depth -2.6pt \kern -6pt \intop}\nolimits_{#1}}
{\mathop{\vrule width 5pt height 3 pt depth -2.6pt \kern -6pt
\intop}\nolimits_{#1}} {\mathop{\vrule width 5pt height 3 pt depth
-2.6pt \kern -6pt \intop}\nolimits_{#1}}}

 \newcommand{\R}{\mathbb{R}}

 \newcommand{\Z}{\mathbb{Z}}

 \newcommand{\wM}{\widetilde{M}}

 \newcommand{\dy}{\text{\rm d}y}

\newcommand{\Var}{\text{\rm Var}}

\begin{document}

\title[Discrete Tanaka's Theorem]{On a discrete version of Tanaka's theorem for maximal functions}

\author[Bober]{Jonathan Bober}
\address{School of Mathematics, Institute for Advanced Study, Einstein Drive, Princeton, NJ, 08540.}
\email{bober@math.ias.edu}

\author[Carneiro]{Emanuel Carneiro}
\address{School of Mathematics, Institute for Advanced Study, Einstein Drive, Princeton, NJ, 08540.}
\email{ecarneiro@math.ias.edu}

\author[Hughes]{Kevin Hughes}
\address{Department of Mathematics, Princeton University, Fine Hall, Washington Road, Princeton, NJ, 08544}
\email{kjhughes@math.princeton.edu}

\author[Pierce]{Lillian B. Pierce}
\address{School of Mathematics, Institute for Advanced Study, Einstein Drive, Princeton, NJ, 08540.}
\email{lbpierce@math.ias.edu}


%    General info

\date{January 11, 2011}


\begin{abstract}
In this paper we prove a discrete version of Tanaka's theorem \cite{Ta} for the Hardy-Littlewood maximal operator in dimension $n=1$, both in the non-centered and centered cases. For the non-centered maximal operator $\wM  $ we prove that, given a function $f: \Z \to \R$ of bounded variation, %L
$$\Var(\wM  f) \leq \Var(f),$$
where $\Var(f)$ represents the total variation of $f$. For the centered maximal operator $M$ we prove that, given a function $f: \Z \to \R$ such that $f \in \ell^1(\Z)$,  %L
$$\Var(Mf) \leq C \|f\|_{\ell^1(\Z)}.$$
This provides a positive solution to a question of Haj\l asz and Onninen \cite{HO} in the discrete one-dimensional case. 
\end{abstract}

\maketitle

\section{Introduction}

It is natural to expect that an averaging operator should have certain smoothing properties; for instance, the spherical means on $\R^d$ map $L^2$ to $W^{\frac{d-1}{2},2}$ \cite[ Chapter 8, \S5.21]{SteinHA}. So one could expect that a maximal operator, being a supremum over averages, should not behave too differently. In fact, if maximal operators are not smoothing operators, at least they do not destroy the regularity of functions, up to one weak derivative. This is the principle behind the program started in 1998 by Kinunnen \cite {Ki} that studies the regularity of maximal operators acting on Sobolev functions. Since then, many authors have contributed to extend the theory, for instance \cite{CM},  \cite{HM}, \cite{KL}, \cite{KiSa}, \cite{Lu1}, always having in the background the general principle that, for maximal operators, an $L^p$-bound implies a $W^{1,p}$-bound.

Things become more difficult when one works with $L^1$-functions, since the Hardy-Littlewood maximal operator does not map $L^1$ to $L^1$. For $f \in L^1_{loc}(\R^n)$ we define the {\it centered} maximal operator as follows:
\begin{equation*}
\mathcal{M}f(x) = \sup_{r>0}\frac{1}{m(B(x,r))} \int_{B(x,r)} |f(y)|\, \dy \,,
\end{equation*}
where $B(x,r)$ is the ball in $\R^n$ centered at $x$ with radius $r$ and $m(B(x,r))$ is the $n$-dimensional Lebesgue measure of this ball. In 2004, Haj\l asz and Onninen \cite[Question 1]{HO} asked the following question: \\
\\
{\bf Question A}: Is the operator $f \mapsto |\nabla \mathcal{M} f|$ bounded from $W^{1,1}(\R^n)$ to $L^1(\R^n)$?\\
\\
Observe that by dilation invariance, a bound of the type
\begin{equation*}
\|\nabla \mathcal{M}f\|_{L^1(\R^n)} \leq C \left( \|f\|_{L^1(\R^n)} + \|\nabla f\|_{L^1(\R^n)} \right)
\end{equation*}
implies that 
\begin{equation}\label{Tan1}
\|\nabla \mathcal{M}f\|_{L^1(\R^n)} \leq C \|\nabla f\|_{L^1(\R^n)}\,,
\end{equation}
and thus the fundamental question here is to compare the variation of the maximal function with the variation of the original function (perhaps having the additional information that $f \in L^1(\R^n))$.

In Tanaka's elegant paper \cite{Ta}, he gave a positive answer to this question for the {\it non-centered} maximal operator in dimension $n=1$. For $f \in L^1_{loc}(\R^n)$ the {\it non-centered} maximal operator is defined as follows:
\begin{equation*}
\widetilde{\mathcal{M}} f(x) = \sup_{\stackrel{r>0}{x \in B_r}}\frac{1}{m(B_r)} \int_{B_r} |f(y)|\, \dy \,,
\end{equation*}
where the supremum is now taken over all balls $B_r$ simply containing $x$.
Tanaka established (\ref{Tan1}) for $\widetilde{\mathcal{M}}$ when $f \in W^{1,1}(\R)$, with constant $C=2$. This result was later refined by Aldaz and P\'{e}rez L\'{a}zaro \cite[Theorem 2.5]{AP} who obtained the sharp constant $C=1$, under only the assumption that $f$ has bounded variation.

Philosophically, the {\it non-centered} version is a smoother operator than the centered version since it contains more averages,  making it easier to handle. Nevertheless, we expect that (\ref{Tan1}) should also hold in the case of the {\it centered} maximal operator with constant $C=1$, when $n=1$ (one may take $f (x)= \chi_{[a,b]}(x)$ for an extremal example). So far, however, Question A remains untouched for the {\it centered} version, even in the case $n=1$.

\subsection{The discrete one-dimensional setting} 
Finding discrete analogues for $L^p$-bounds in harmonic analysis is a topic of ongoing research. In the simplest cases, $\ell^p$-bounds for discrete analogues of classical operators such as Calder\'{o}n-Zygmund singular integral operators, fractional integral operators, and the Hardy-Littlewood maximal function follow from known $L^p$-bounds for the original operators in the Euclidean setting,  via elementary comparison arguments (see \cite{SW1}, \cite{SW2}). But $\ell^p$-bounds for discrete analogues of more complicated operators, such as singular, fractional, and maximal Radon transforms (involving integration over a submanifold, or family of submanifolds), are not implied by results in the continuous setting, and moreover the discrete analogues are resistant to conventional methods. Indeed, discrete operators may even behave differently from their continuous counterparts, as is exhibited by the discrete spherical maximal operator \cite{MSW}.
It is only recently that substantial progress has been made on discrete operators with Radon characteristics via techniques motivated by the circle method of Hardy and Littlewood, a technique from number theory  pioneered in the context of discrete analogues by Bourgain \cite{Bour88A}, \cite{Bour88B}, and  further developed in a number of interesting cases (see for example \cite{IMSW}, \cite{IW}, \cite{MSW}, \cite{Pie10},  \cite{SW1}, \cite{SW2}).

In this paper we introduce the study of the regularity theory of discrete maximal operators in one dimension.
Let $f: \Z \to \R$ be a discrete function and let $\Z^+ = \{0, 1, 2, 3, \ldots,\}$.
The discrete {\it centered} Hardy-Littlewood maximal operator is defined by
\begin{equation*}
 Mf(n) = \sup_{r \in \Z^{+}} \frac{1}{(2r + 1)} \sum_{k=-r}^{k=r} |f(n+k)|\,,
\end{equation*}
while the {\it non-centered} version is defined by
\begin{equation*}
 \wM  f(n) = \sup_{r,s \in \Z^{+}} \frac{1}{(r +s+ 1)} \sum_{k=-r}^{k=s} |f(n+k)|\,.
\end{equation*}
Our aim is to answer discrete analogues of Question A for these operators. They clearly do not belong to the Radon transform paradigm, and we will not call upon the circle method; instead the challenge arises, at least in the case of the centered maximal operator $M$, from the fact that the analogous result in the continuous setting is not yet known!

In order to study regularity properties of discrete operators, we establish the following conventions.
For $1 \leq p <\infty$, the $\ell^p$-norm of a function $f: \Z \to \R$ is
\begin{equation*}
 \|f\|_{\ell^p(\Z)} = \left(\sum_{n=-\infty}^{\infty} |f(n)|^p \right)^{1/p}\,,
\end{equation*}
and the $\ell^{\infty}$-norm is 
\begin{equation*}
 \|f\|_{\ell^{\infty}(\Z)} = \sup_{n \in \Z} |f(n)|.
\end{equation*}
We define the derivatives of a discrete function by 
\begin{align*}
 f'(n) &= f(n+1) - f(n)\,,\\
f''(n) & = f(n+2) - 2f(n+1) + f(n)\,,\\
f'''(n) & = f(n+3) - 3f(n+2) + 3 f(n+1) - f(n)\,,
\end{align*}
and so on. 
The space corresponding to $W^{k,p}(\R)$ is then defined to be the set of discrete functions with finite $w^{k,p}(\Z)$-norm, where
\[ ||f||_{w^{k,p}(\Z)}  = \sum_{j=0}^k ||f^{(j)}||_{\ell^p(\Z)}.\]
But note that by the triangle inequality, for any $k \geq 1$, 
\begin{equation}\label{TI}
 \|f^{(k)}\|_{\ell^p(\Z)} \leq 2^k \|f\|_{\ell^p(\Z)} \,,
\end{equation}
thus in the discrete setting, any $\ell^p$-bound automatically provides a $w^{k,p}$-bound, for any $k \geq 1$ (in fact, the discrete $w^{k,p}$-spaces are just the classical $\ell^p$-spaces with an equivalent norm).  
This might make our efforts to transfer the regularity theory for maximal operators to the discrete setting seem almost vacuous.

However, the situation is highly nontrivial when we deal with $\ell^1$-functions or functions of bounded variation.
We define the total variation of $f:\Z \to \R$ by
\begin{equation*}
\Var(f) = \|f'\|_{\ell^1(\Z)} = \sum_{n=-\infty}^{\infty} |f(n+1) - f(n)|\,.
\end{equation*}
Our first result, a discrete version of Tanaka's theorem for the discrete {\it non-centered} maximal operator, with sharp constant, is as follows:

\begin{theorem}\label{thm1}
Let $f: \Z \to \R$ be a function of bounded variation. Then
\begin{equation*}
\Var (\wM  f) \leq \Var(f)\,,
\end{equation*}
and the constant $C=1$ is the best possible.
\end{theorem}
We shall prove this result in Section 2, adapting some of the ideas of the original proof of Tanaka for the continuous case. 
It is not hard to see that the constant $C=1$ is best possible in Theorem \ref{thm1}, for it suffices to consider the function
\begin{equation}\label{ex1}
f(n) = \left\{
\begin{array}{cc}
1& \textrm{if} \ n=0,\\
0& \textrm{otherwise}.
\end{array}
\right.
\end{equation}

Dealing with the {\it centered} maximal operator is a much more subtle and intricate problem. 
By an extensive analysis of examples we are led to believe that the same bound should hold for the {\it centered} maximal operator:\\
\\
{\bf Question B}: Let $f: \Z \to \R$ be a function of bounded variation. Is it true that %L
\begin{equation}\label{QB}
\Var(Mf) \leq \Var(f)\ ?\\
\end{equation}
Motivated by Question B, we prove the following result.
 
\begin{theorem}\label{thm2}
Let $f: \Z \to \R$ be a function in $\ell^1(\Z)$. Then 
\begin{equation}\label{Var1}
 \Var (Mf) \leq \left(2 + \frac{146}{315} \right) \|f\|_{\ell^1(\Z)}.
\end{equation}
\end{theorem}

Theorem \ref{thm2} represents partial progress toward Question B. In fact, from (\ref{TI}), inequality (\ref{QB}) would imply (\ref{Var1}) with constant $C=2$, which would be sharp, with an extremal example given by (\ref{ex1}). The constant we obtain here is slightly bigger 
$$C = \left(2 + \frac{146}{315} \right)  = 2 \left(1 + \frac{1}{5} + \frac{1}{7} - \frac{1}{9}\right),$$
and it appears due to combinatorial arguments in our proof (see Lemma \ref{lemma2} below). It is an interesting question whether one can improve this constant towards the conjectured sharp value.

We expect higher dimensional analogues of these results, both in the continuous and discrete cases, to hold as well (see the original question by Haj\l asz and Onninen \cite{HO}). However, the methods of Tanaka \cite{Ta} and Aldaz and P\'{e}rez L\'{a}zaro \cite{AP} for the one dimensional continuous (uncentered) case, and ours for the discrete (centered and uncentered) cases, do not easily adapt to higher dimensions.

The simplicity and innocence of the objects and statements described above might appear misleading at first glance. Before moving to the proofs, we encourage the interested reader to familiarize her/himself with the discrete maximal problem, especially Question B above, in order to better appreciate the beauty and the difficulties of the interplay between analysis and combinatorics, still not completely understood, in this problem. 

\section{Proof of Theorem \ref{thm1}}
Since $\Var(|f|) \leq \Var(f)$ we may assume without loss of generality that $f$ takes only non-negative values. A function of bounded variation will certainly be bounded and thus, at each point $n$, the averages will also be bounded. However, since we do not assume $f  \in \ell^1(\Z)$, we must be aware of the fact that the supremum over these averages might not be realized.

We will say that a point $n$ is a {\it local maximum} of $f$ if
\begin{equation*}
f(n-1) \leq f(n)  \ \ \ \textrm{and} \ \ \ f(n) > f(n+1).
\end{equation*}
Similarly, a point $n$ is a {\it local minimum} of $f$ if
\begin{equation*}
f(n-1) \geq f(n) \ \ \ \textrm{and} \ \ \ f(n) < f(n+1).
\end{equation*}
The following lemma identifies a key property of the local maxima of $\wM  f$. 

\begin{lemma}\label{lemma_extrema}
If $n$ is a local maximum of $\wM  f$, then $\wM  f(n) = f(n)$.
\end{lemma}
\begin{proof}
Let $n \in \Z$ be a point such that $\wM  f(n) > f(n)$ and assume that $n$ is a local maximum of $\wM  f$. Let us arrive at a contradiction.

\medskip

\noindent{\it Case 1.} {\it $\wM  f(n)$ is attained for some $r,s \in \Z^+$}. Let $I = [n-r,n+s]$ be such an interval. If $s\geq 1$, then we may take an average over the same interval $I$ in order to bound $\wM  f(n+1)$, i.e.
$$ \wM  f(n+1) \geq \frac{1}{(r +s+ 1)} \sum_{k=-r}^{k=s} f(n+k) = \wM  f(n),$$
a contradiction. If $s=0$, since $\wM  f(n) > f(n)$, we must have $r \geq 1$. We then have
\begin{align*}
\wM  f(n)  & = \frac{1}{(r + 1)} \sum_{k=-r}^{k=0} f(n+k),\\
\wM  f(n-1) & \geq \frac{1}{r} \sum_{k=-r}^{k=-1} f(n+k).
\end{align*}
This also leads to a contradiction by observing that 
$$\wM  f(n) \leq (r+1) \wM  f(n) - r \wM  f(n-1)  \leq f(n). $$

\medskip

\noindent{\it Case 2.} {\it $\wM  f(n)$ is not attained for any $r,s \in \Z^+$}. In this case we see that $\wM  f(m) \geq \wM  f(n)$ for any $m \in \Z$, which is a contradiction. In fact, given $\varepsilon >0$, we may take an interval $I = [n-r,n+s]$ big enough such that 
$$\wM  f(n) -\varepsilon \leq \frac{1}{(r +s+ 1)} \sum_{k=-r}^{k=s} f(n+k).$$
For a fixed $m \in \Z$, we consider the average over the interval $I' = [m-r,m+s]$. Letting $C = \|f\|_{\ell^{\infty}(\Z)}$, we obtain
\begin{align*}
\wM  f(m) & \geq \frac{1}{(r +s+ 1)} \sum_{k=-r}^{k=s} f(m+k)\\
& \geq \big(\wM  f(n) -\varepsilon\big) - \frac{2C |m-n|}{r + s +1}.
\end{align*}
Letting $(r + s) \to \infty$ we get
$$\wM  f(m) \geq (\wM  f(n) -\varepsilon\big).$$
Letting $\varepsilon \to 0$, we arrive at the desired conclusion.
\end{proof}

We now finish the proof of Theorem \ref{thm1}. From now on let  us consider the alternating sequence of local maxima $\{a_i\}_{i\in \Z}$ and local minima $\{b_i\}_{i\in \Z}$ of $\wM  f$, satisfying
\begin{equation}\label{sequence0}
...< b_{-2} < a_{-2} < b_{-1} < a_{-1} < b_0 < a_0 < b_1 < a_1 < b_2 < a_ 2 < ....
\end{equation}
The sequence (\ref{sequence0}) can be finite or infinite depending on the behavior of the tails of $\wM  f$. Let us consider the different cases.\\
\\
{\it Case 1}. The sequence (\ref{sequence0}) is infinite on both sides.\\

In this case we have 
\begin{align}\label{VarNC}
\begin{split}
\Var\big(\widetilde{M}f\big) &= \sum_{i=-\infty}^{\infty} \Big\{ \widetilde{M}f(a_{i-1}) - \widetilde{M}f(b_{i})\Big\} +  \Big\{ \widetilde{M}f(a_{i}) - \widetilde{M}f(b_{i})\Big\}\\
& = \sum_{i=-\infty}^{\infty} \Big\{ f(a_{i-1}) - \widetilde{M}f(b_{i})\Big\} +  \Big\{f(a_{i}) - \widetilde{M}f(b_{i})\Big\}\\
& \leq \sum_{i=-\infty}^{\infty} \big\{ f(a_{i-1}) - f(b_{i})\big\} +  \big\{f(a_{i}) - f(b_{i})\big\}  \leq \Var(f).\\
\end{split}
\end{align}
{\it Case 2}. The sequence (\ref{sequence0}) is finite on one (or both) side(s).\\

In this case several different behaviors might occur, but they are essentially treated in the same way, using (\ref{VarNC}) and a minor modification in the tail(s). Suppose for instance that $a_k$ is the last local maximum. The function $\wM  f(n)$ must be monotonically non-increasing for $n \geq a_k$ and since it is bounded, the limit
\begin{equation*}
\wM  f(\infty) = \lim_{n \to \infty} \wM  f(n) = c
\end{equation*}
will exist. In this case we will have
\begin{equation*}
\liminf_{n\to \infty}f(n) \leq c.
\end{equation*}
Below we write
\begin{equation*}
 \Var( f)_{[a,b]} = \sum_{n=a}^{b-1} |f(n+1) - f(n)|\,,
\end{equation*}
for the variation of $ f$ on the interval $[a,b]$, where $a$ and $b$ are integers (or possibly $\pm \infty$). Therefore we have
\begin{align*}
\begin{split}
&\Var\big(\widetilde{M}f\big) = \Var\big(\widetilde{M}f\big)_{[-\infty, a_k]} + \Var\big(\widetilde{M}f\big)_{[a_k,\infty]} \\
& = \sum_{i = -\infty}^{k} \left\{ \Big(\widetilde{M}f(a_{i-1}) - \widetilde{M}f(b_{i})\Big) +  \Big( \widetilde{M}f(a_{i}) - \widetilde{M}f(b_{i}) \Big) \right\} + \Big(\widetilde{M}f(a_k) - c \Big)\\
& = \sum_{i = -\infty}^{k} \left\{ \Big(f(a_{i-1}) - \widetilde{M}f(b_{i})\Big) +  \Big( f(a_{i}) - \widetilde{M}f(b_{i}) \Big) \right\} + \big(f(a_k) - c \big)\\
& \leq  \sum_{i = -\infty}^{k} \left\{ \big(f(a_{i-1}) - f(b_{i})\big) +  \big( f(a_{i}) - f(b_{i}) \big) \right\} + \big(f(a_k) - c \big)\\
& \leq \Var(f)_{[-\infty, a_k]} + \Var(f)_{[a_k,\infty]} = \Var(f).
\end{split}
\end{align*}
The argument for all the other cases is a minor modification of this one. This concludes the proof of Theorem \ref{thm1}.

\section{Proof of Theorem \ref{thm2}}
One can begin consideration of the discrete {\it centered} maximal operator by investigating whether Lemma \ref{lemma_extrema}, or any natural modification of it, continues to hold. The following example shows that this need not be the case: 
\begin{equation*}\label{ex2}
f(n) = \left\{
\begin{array}{cc}
10& \textrm{if} \ n=\pm4,\\
0& \textrm{otherwise}.
\end{array}
\right.
\end{equation*}
One should not expect the local maxima of $Mf$ to touch $f$, or even expect that $Mf$ should be convex in each interval in which it disconnects from $f$. 

Thus new ideas are required to approach this problem. We start again by assuming that $f$  takes only non-negative values, and consider the sequence of local maxima $\{a_i\}_{i\in \Z}$ and local minima $\{b_i\}_{i\in \Z}$ of $Mf$, satisfying
\begin{equation}\label{sequence}
...< b_{-2} < a_{-2} < b_{-1} < a_{-1} < b_0 < a_0 < b_1 < a_1 < b_2 < a_ 2 < ....
\end{equation}
We have 
\begin{equation}\label{Var}
 \Var(Mf) = 2\sum_{i=-\infty}^{\infty} \left(Mf(a_i) - Mf(b_i)\right).
\end{equation}
{\bf Remark}: If the sequence (\ref{sequence}) terminates on one or both ends, we modify the sum (\ref{Var}) accordingly as follows. Since $f \in \ell^1(\Z)$ we must have $\lim_{n \to \pm \infty} Mf(n) = 0$, and this implies that if the sequence terminates, it would terminate with a last maximum $a_k$ and/or a first maximum $a_l$ (i.e. it would not terminate with a minimum). If there is a first maximum $a_l$ we consider 
\begin{equation*}
 \Var(Mf) = 2 Mf(a_l) + 2 \sum_{i=l+1}^{\infty} \left(Mf(a_i) - Mf(b_i)\right),
\end{equation*}
and make minor modifications in the argument below; similar modifications apply if there is a last maximum $a_k$.\\

For each local maximum $a_i$ we let $r_i$ be the smallest radius such that 
\begin{equation}\label{avmax}
Mf(a_i) = A_{r_i}f(a_i) = \frac{1}{(2r_i + 1)} \sum_{k=-r_i}^{k=r_i} f(a_i+k)\,,
\end{equation}
where we denote by $A_{r}$ the averaging operator of radius $r$ (since $f \in \ell^1(\Z)$ this radius exists). For each point $b_i$ we consider the average of radius $s_i = r_i + (a_i - b_i)$ and since we have
\begin{equation}\label{avmin}
 Mf(b_i) \geq A_{s_i}f(b_i)\,,
\end{equation}
it follows from (\ref{Var}), (\ref{avmax}) and (\ref{avmin}) that 
\begin{equation}\label{VarSum}
  \Var(Mf) \leq 2\sum_{i=-\infty}^{\infty} \left(A_{r_i}f(a_i) - A_{s_i}f(b_i)\right).
\end{equation}
Observe that the interval $[b_i - s_i, b_i + s_i]$ contains the interval $[a_i - r_i, a_i + r_i]$ and they both have the same right endpoint. Now we fix an integer $n$ and we will evaluate the maximum contribution that $f(n)$ can give to the sum on the right hand side of (\ref{VarSum}). 

For each $i \in \Z$ , if $n \in [a_i - r_i, a_i + r_i]$, then $n \in [b_i - s_i, b_i + s_i]$ and $f(n)$ contributes to  $\left(A_{r_i}f(a_i) - A_{s_i}f(b_i)\right)$ the amount %L
\begin{equation}\label{contribution of n}
\frac{f(n)}{2r_i + 1} - \frac{f(n)}{2( r_i + (a_i - b_i)) + 1}.
\end{equation}
If $n \notin [a_i - r_i, a_i + r_i]$ the contribution of $f(n)$ to $\left(A_{r_i}f(a_i) - A_{s_i}f(b_i)\right)$ is zero or even negative and we disregard it. Now observe that if the contribution (\ref{contribution of n}) occurs we must have $r_i \geq |n - a_i|$, and therefore one can show that
\begin{align}\label{upperbound}
\begin{split}
f(n) &\left(\frac{1}{2r_i + 1} - \frac{1}{2( r_i + (a_i - b_i)) + 1}\right)\\
& \ \ \ \ \ \ \   \leq f(n) \left(\frac{1}{2|n-a_i| + 1} - \frac{1}{2( |n-a_i| + (a_i - b_i)) + 1}\right)\\
& \ \ \ \ \ \ \  \ \ \ \ \ \ \ \ \   \leq  f(n) \left(\frac{1}{2|n-a_i| + 1} - \frac{1}{2( |n-a_i| + (a_i - a_{i-1})) + 1}\right),
\end{split}
\end{align}
where in the last step we just used the ordering (\ref{sequence}). If we sum (\ref{upperbound}) over the index $i$ we obtain an upper bound for the total contribution of $f(n)$ to the right hand side of (\ref{VarSum}), namely
\begin{equation}\label{FinalSum}
 2 f(n) \sum_{i=-\infty}^{\infty}\left(\frac{1}{2|n-a_i| + 1} - \frac{1}{2( |n-a_i| + (a_i - a_{i-1})) + 1}\right).
\end{equation}
Theorem \ref{thm2} will follow if we prove that for {\it any} strictly increasing sequence $\{a_i\}_{i \in \Z}$ of integers, the sum in (\ref{FinalSum}) is bounded by a universal constant $C$. This is proved in Lemma \ref{lemma2} below. To conclude the proof of Theorem \ref{thm2}, we will ultimately sum the maximum contributions of all $f(n)$'s to the total variation (\ref{VarSum}) of $Mf$ to prove, as desired, that
\begin{equation}\label{eq-f}
 \Var(Mf) \leq 2 C \sum_{n=-\infty}^{\infty} f(n) = 2 C\,\|f\|_{\ell^1(\Z)}.
\end{equation}

\begin{lemma}\label{lemma2}
Given $n\in \Z$, for {\it any} strictly increasing sequence $\{a_i\}_{i \in \Z}$ of integers,
\begin{equation*}
 \sum_{i=-\infty}^{\infty}\left(\frac{1}{2|n-a_i| + 1} - \frac{1}{2( |n-a_i| + (a_i - a_{i-1})) + 1}\right)\leq \frac{4}{3}. 
 \end{equation*}
If furthermore $a_i - a_{i-1} \geq 2$ for all $i\in \Z$, the constant $C = \frac{4}{3}$ may be replaced by $C = 1 + \frac{1}{5} + \frac{1}{7} - \frac{1}{9}.$
\end{lemma}
\begin{proof}
It is sufficient to prove the result for $n=0$ since we can shift any sequence $a_i \mapsto a_i +n$. For $n=0$ we aim to prove that 
\begin{equation}\label{nicesum}
S = \sum_{i=-\infty}^{\infty}\left(\frac{1}{2|a_i| + 1} - \frac{1}{2( |a_i| + (a_i - a_{i-1})) + 1}\right)\leq C\,.
\end{equation}
By shifting the indices we can also assume that $a_{-1} \leq 0 < a_0$. We divide our sum (\ref{nicesum}) into two parts
\begin{align*}
\begin{split}
S &=  \sum_{i=-\infty}^{-1}\left(\frac{1}{-2a_i + 1} - \frac{1}{-2a_{i-1} + 1}\right) + \sum_{i=0}^{\infty}\left(\frac{1}{2a_i + 1} - \frac{1}{2( a_i + (a_i - a_{i-1})) + 1}\right)\\
\\
& = S_1 + S_2.
\end{split}
\end{align*}
The first sum $S_1$ is a telescoping sum and we find that
\begin{equation*}
 S_1 \leq \frac{1}{-2a_{-1} + 1}.
\end{equation*}
(This continues to hold if the sequence terminates to the left as $i \to -\infty$.)

The second sum is more involved and we use the following inequality, for  integers $m > n\geq 0$:
\begin{equation}\label{ineqtricky}
 \frac{1}{2m + 1} - \frac{1}{2( m + (m - n)) + 1} \leq \frac{1}{2(n+1) +1} - \frac{1}{2(m+1) + 1}.
\end{equation}
Inequality (\ref{ineqtricky}) can be proved simply by clearing denominators and observing that $m \geq n+1$.
We then use (\ref{ineqtricky}) to bound $S_2$ as follows:
\begin{align*}
S_2 &=  \left(\frac{1}{2a_0 + 1} - \frac{1}{2( a_0 + (a_0 - a_{-1})) + 1}\right) + \sum_{i=1}^{\infty}\left(\frac{1}{2a_i + 1} - \frac{1}{2( a_i + (a_i - a_{i-1})) + 1}\right)\\
\\
& \leq \left(\frac{1}{2a_0 + 1} - \frac{1}{2( a_0 + (a_0 - a_{-1})) + 1}\right) + \sum_{i=1}^{\infty} \left(\frac{1}{2(a_{i-1}+1) +1} - \frac{1}{2(a_i+1) + 1}\right)\\
\\
& \leq  \left(\frac{1}{2a_0 + 1} - \frac{1}{2( a_0 + (a_0 - a_{-1})) + 1}\right)+ \frac{1}{2(a_{0}+1) +1}.\\
\end{align*}
(This also continues to hold if the sequence terminates to the right as $i \to \infty$.)
We have thus arrived at 
\begin{align}\label{Final1}
\begin{split}
S & = S_1 + S_2 \\
& \leq  \frac{1}{-2a_{-1} + 1} + \left(\frac{1}{2a_0 + 1} - \frac{1}{2( a_0 + (a_0 - a_{-1})) + 1}\right)+ \frac{1}{2(a_{0}+1) +1}.
\end{split}
\end{align}
Recall that here $a_{-1} \leq 0 < a_0$ are integers. For any $0<a_0$ it is easy to see that \eqref{Final1} is maximized when $a_{-1}=0$. Then a simple analysis of cases yields that $a_0=1$ is the maximal choice, proving that \eqref{Final1} is bounded by 
$$C = 1 + \frac{1}{3} = \frac{4}{3}.$$
If we impose the condition $a_0 - a_{-1} \geq 2$, again we can easily see that the maximum of (\ref{Final1}) occurs when $a_{-1} = -1$ or $0$. An analysis of a few cases confirms that \eqref{Final1} is maximized when $a_{-1} = 0$  and $a_0 = 2$, giving the upper bound 
\begin{equation}\label{C2}
C = 1 + \frac{1}{5} + \frac{1}{7} - \frac{1}{9}= 1 + \frac{73}{315}\,,
\end{equation}
and this finishes the proof of the lemma.
\end{proof}

Observe that in the setting of our maximal operator, we have the condition $a_i - a_{i-1} \geq 2$ for all $i\in \Z$ (there must be a local minimum between any two consecutive local maxima), and thus we can use Lemma \ref{lemma2} with constant  $C$ given by \eqref{C2} in expression \eqref{eq-f} to conclude that 
\begin{equation*}
 \Var(Mf) \leq 2 \left(1 + \frac{73}{315}\right) \,\|f\|_{\ell^1(\Z)}.
\end{equation*}

\section*{Acknowledgments}
The authors would like to thank Jeffrey Vaaler, Jean Bourgain, Diego Moreira and Dimitris Koukoulopoulos for helpful comments during the preparation of this work. J. Bober, E. Carneiro and L. B. Pierce  acknowledge support from the Institute for Advanced Study and the National Science Foundation under agreement No. DMS-0635607. 
%Any opinions, findings and conclusions or recommendations expressed in this material are those of the authors and do not necessarily reflect the views of the National Science Foundation. 
E. Carneiro also acknowledges support from CAPES/FULBRIGHT grant BEX 1710-04-4.
L. B. Pierce is also funded by the Simonyi Fund and National Science Foundation grant DMS-0902658.

\section*{Remarks}

Our previous version of this manuscript (arXiv:1005.3030v3), which happens to be the published version [Proc. Amer. Math. Soc. 140 (2012), 1669-1680], had an oversight in the proof of Lemma \ref{lemma_extrema} (iii), when we claimed that 
\begin{equation}\label{eq_lateral}
\wM  f(n) = \max\{M_Lf(n), M_Rf(n)\}, \ \ \forall n \in \Z,
\end{equation}
which is not true. The previous proof of Lemma \ref{lemma_extrema} (iii), and hence of Theorem \ref{thm1}, was correct for the slightly different discrete non-centered maximal operator
\begin{equation*}
 \widetilde{\wM}  f(n) = \sup_{r,s \in \Z^{+}} \frac{1}{(r +s+ \lambda_{rs})}\left\{ \left( \sum_{k=-r}^{k=-1} |f(n+k)|\right) + \lambda_{rs} |f(n)| + \left(\sum_{k=1}^{k=s} |f(n+k)|\right)\right\}\,,
\end{equation*}
where $\lambda_{rs} = 1$ if $rs>0$ and $\lambda_{rs} = \frac12$ if $rs=0$ (i.e. when we take an average over an interval that has $n$ as an endpoint, we consider is contribution with weight $\frac12$). We have now briefly adjusted our proof to address our original claim for the operator $\wM$. In fact, Lemma \ref{lemma_extrema} (and Theorem \ref{thm1}) hold for other variants of non-centered discrete maximal operators (e.g. placing weight $1/2$ in one or both endpoints of the interval, considering only intervals with an odd number of integer points, considering lateral operators, etc.). The strategy of the proof is the same, with minor adjustments. The oversight in Lemma  \ref{lemma_extrema} (iii) in our previous version was observed independently by Dariusz Kosz and Matt Rosenzweig, to whom we are thankful.

\bibliographystyle{amsplain}

\end{document}